# МЕТОД ДИСКРЕТИЗАЦИИ ДРОБНО-ПРОИЗВОДНЫХ ЛИНЕЙНЫХ СИСТЕМ ОБЫКНОВЕННЫХ ДИФФЕРЕНЦИАЛЬНЫХ УРАВНЕНИЙ С ПОСТОЯННЫМИ КОЭФФИЦИЕНТАМИ

**Алиев Фикрет А., Алиев Н.А., Велиева Н.И., Гасымова К.Г., Мамедова Е.В.**

**Абстракт.** Развивается точный метод дискретизации для решения линейных систем обыкновенных дробно-производных дифференциальных уравнений с постоянными матричными коэффициентами (ЛСОДПДУПК). Показывается, что полученная линейная дискретная система в данном случае не имеет постоянные матричные коэффициенты. Далее данный метод сравнивается с известным приближенным методом. Приведенная схема развивается для произвольных линейных систем с кусочно-постоянными возмущениями. Результаты применяются для дискретизации линейных управляемых систем и иллюстрируются числовыми примерами.

**Введение**

Как известно [1-3], задача Коши для систем линейных дифференциальных уравнений с постоянными матричными коэффициентами

$$\dot{x} = Ax, \quad x(t_0) = x_0, \qquad (1)$$

имеет решение в следующем виде

$$x(t) = e^{A(t-t_0)}x(t_0), \qquad (2)$$

где $A$-квадратная постоянная матрица размерности $n \times n$, $x - n$-мерный неизвестный вектор.

Соотношение (2) позволяет вычислить численно решение $x(t)$ в каждой точке $t_i$, входящей в область определения $[t_0, T]$. Действительно, если разделить область определения $[t_0, T]$ с постоянным шагом $\Delta$, то из (2) для каждой точки $t_i$ и $t_{i+1}$, имеем следующее разностное соотношение

$$x(t_{i+1}) = e^{A\Delta}x(t_i), \quad x(t_0) = x_0, \qquad (3)$$

Если удастся точно вычислить матрицу $e^{A\Delta}$, то разностное соотношение (3) на каждой точке $t_i$ определяет $x(t_i)$ точно. Если шаг $\Delta$ не равномерно (не постоянно) и имеет вид $\Delta_i$, то (3) заменяется следующим разностным уравнением [1-3, 4]

$$x(t_{i+1}) = e^{A\Delta_i} x(t_i), \qquad (4)$$

которое является системой нестационарных линейных конечно-разностных уравнений.

Если рассмотреть дискретизацию задачи Коши (1) заменив производные $\dot{x}$ разностным выражением $\dfrac{x(t_{i+1}) - x(t_i)}{\Delta}$, то приближенная форма уравнения (1) будет иметь вид

$$\tilde{x}(t_{i+1}) = (E + \Delta \cdot A)\tilde{x}(t_i), \quad \tilde{x}(t_0) = x_0, \qquad (5)$$

где при оптимальном выборе $\Delta$ можно обеспечить достаточную близость решения (5) к (3). Легко видеть, что (5) фактически является линеаризацией матриц $e^{A\Delta}$ из соотношения (3). Отметим, что как отмечено в [5], надо выбрать $\Delta$ так, чтобы структурное свойство уравнения (1) не менялось при переходе к (5).

## 2. Приближенная дискретизация для дробных производных через фундаментальные матрицы

Если рассмотреть системы обыкновенных дробно-производных дифференциальных уравнений с постоянными матричными коэффициентами вида

$$D^{\alpha} x(t) = Ax(t), \quad x(t_0) = x_0, \qquad (6)$$

в интервале $[t_0, T]$, то картина дискретизации (4) для данного случая поменяется, т.е. простая формула (4) (или же (5)) для случая (6) полностью изменяется. Действительно, как известно [6] (см. формулы (3.109) при $B = 0$), приближенное решение уравнения (6) в интервале определения $[t_0, T]$ имеет вид

$$x((k+1)\Delta) = (A\Delta^{\alpha} + \alpha E) x(k\Delta) - \sum_{i=2}^{k+1} (-1)^i \binom{\alpha}{i} x((k+1-i)\Delta), \qquad (7)$$

и для $k = 0$

$$x(\Delta) = (A\Delta^{\alpha} + \alpha E) x(0), \qquad (8)$$

где

$$\binom{\alpha}{i} = \frac{\alpha(\alpha-1)(\alpha-2)\ldots(\alpha-i+1)}{i!}$$

Используя результаты [6] приведем приближенное вычисление между $x(k+1)$ и $x(k)$. Как известно [6] между $x(k+1)$ и $x(0)$ имеется следующее соотношение

$$x(k+1) = \Phi(k+1)x(0) \qquad (9)$$

где

$$\Phi(k+1) = (A\Delta^\alpha + \alpha E)\Phi(k) - \sum_{i=2}^{k+1}(-1)^i \binom{\alpha}{i}\Phi(k+1-i),$$

(10)

$$\Phi(1) = (A\Delta^\alpha + \alpha E)\Phi(0), \quad \Phi(0) = I.$$

Из (9) для $x(0)$ имеем

$$x(0) = \Phi^{-1}(k)x(k) \qquad (10')$$

Как следует из (9), $(10')$ заменяя $k+1$ через k и m соответственно получим

$$x(k) = \Phi(k)x(0)$$

$$x(0) = \Phi^{-1}(m)x(m)$$

получим

$$x(k) = \Phi(k)\Phi^{-1}(m)x(m).$$

Обозначим

$$\Phi(k, m) = \Phi(k) \cdot \Phi^{-1}(m) \qquad (10'')$$

Тогда из $(10'')$ $\Phi(k, m)$ можно принимать как фундаментальная матрица уравнений (6), т.е.

$$x(k) = \Phi(k,m)x(m) \qquad (11)$$

Таким образом, из (11) между $x(k+1)$ и $x(k)$ получим

$$x(k+1) = \Phi(k+1,k)x(k), \qquad (12)$$

которое является приближенной дискретизацией уравнения (6) через фундаментальную матрицу $(10'')$.

## 3. Определение аналитических видов фундаментальной матрицы уравнения (6).

Как показано [7-9], точное решение уравнения (6) (аналогично соотношению (2)) имеет вид

$$x(t) = \left[\sum_{s=0}^{2q} A^{\frac{s}{2p+1}} \frac{t_0^{\frac{s-2q}{2q+1}}}{\frac{s-2q}{2q+1}!}\right]^{-1} \left[\sum_{s=0}^{2q} A^{\frac{s}{2p+1}} \frac{t^{\frac{s-2q}{2q+1}}}{\frac{s-2q}{2q+1}!} + A^{\frac{s+2q+1}{2p+1}} \int_{t_0}^{t} \frac{(t-\tau)^{\frac{s-2q}{2q+1}}}{\frac{s-2q}{2q+1}!} e^{A^{\frac{2q+1}{2p+1}}(\tau-t_0)} d\tau\right] x_0, \qquad (13)$$

где $\alpha = \dfrac{2p+1}{2q+1}$, и $p$, $q$ - натуральные числа.

Таким образом, аналогично дискретное соотношение (3) для задачи (6) можно построить соответствующее соотношение. Тогда из (13) аналитическая форма фундаментальной матрицы для уравнения (6) имеет вид

$$\Phi(t,t_0) = \left[\sum_{s=0}^{2q} A^{\frac{s}{2p+1}} \frac{t_0^{\frac{s-2q}{2q+1}}}{\frac{s-2q}{2q+1}!}\right]^{-1} \left[\sum_{s=0}^{2q} A^{\frac{s}{2p+1}} \frac{t^{\frac{s-2q}{2q+1}}}{\frac{s-2q}{2q+1}!} + A^{\frac{s+2q+1}{2p+1}} \int_{t_0}^{t} \frac{(t-\tau)^{\frac{s-2q}{2q+1}}}{\frac{s-2q}{2q+1}!} e^{A^{\frac{2q+1}{2p+1}}(\tau-t_0)} d\tau\right],$$
$$\Phi(t_0,t_0) = E \qquad (14)$$

Отметим, что (13), (14) позволяет представить общее решение уравнения (6) в виде

$$x(t) = \Phi(t,t_0)x(t_0), \qquad (15)$$

который позволит уже дискретизировать (6) в соответствующем виде (точный аналог (11)).

Далее рассматривается дискретизация уравнения (6) используя соотношения (14), (15), приводится приближенный алгоритм для вычисления $\Phi(t,t_0)$ на основе [6] и сравниваются эти результаты на численном примере. Показывается, что при $\Delta \to 0$ оба результата становятся близкими друг к другу.

Далее результаты развиваются для случая, когда уравнение (6) заменяется уравнением с возмущением

$$D^\alpha x = Ax + f, \quad x(t_0) = x_0, \qquad (16)$$

где $f - n-$ мерный известный вектор возмущений. Используя (13) приводится дискретизация уравнения (16). Отметим, что такая задача более подходящая для линейной дискретной системы управлений [10,15-16], которая требует отдельной разработки.

## 4. Точная дискретизация уравнения (6)

Пусть уравнение (6) определено в интервале $(t_0, T)$. Разделим этот интервал на $\kappa$ частей $(t_i, t_{i+1},)$, $i = 0,1,2,...,k-1$. Тогда на каждом из этих под интервалов дискретный аналог уравнения (6) на основе соотношений (14), (15) имеет вид

$$x(t_{i+1}) = \Phi(t_{i+1}, t_i) x(t_i), \quad x(t_0) = x_0, \qquad (17)$$

$$\Phi(t_{i+1}, t_i) = \left[\sum_{s=0}^{2q} A^{\frac{s}{2p+1}} \frac{t_i^{\frac{s-2q}{2q+1}}}{\left(\frac{s-2q}{2q+1}\right)!}\right]^{-1} \sum_{s=0}^{2q} \left[A^{\frac{s}{2p+1}} \frac{t_{i+1}^{\frac{s-2q}{2q+1}}}{\frac{s-2q}{2q+1}!} + A^{\frac{s+2q+1}{2p+1}} \int_{t_i}^{t_{i+1}} \frac{(t_{i+1}-\tau)^{\frac{s-2q}{2q+1}}}{\frac{s-2q}{2q+1}!} e^{(\tau-t_i)A^{\frac{2q+1}{2p+1}}} d\tau\right] \qquad (18)$$

Отметим, что решений уравнений (6) по формулы (17), (18) дают точные значения в точках $t_i$ $(i = 0,1,...,k)$. Фактически (16), (17) при $\alpha = 1$ (при $p=0$, $q=0$) совпадает с (4).

Вычисление фундаментальных матриц (18) составляет определенные трудности из за входящего в него интеграла. Поэтому подставляя выражение экспоненциальной функции входящего в подынтегральное выражение, как в [9] имеем

$$\Phi(t_{i+1}, t_i) = \left[\sum_{s=0}^{2q} A^{\frac{s}{2p+1}} \frac{t_i^{\frac{s-2q}{2q+1}}}{\left(\frac{s-2q}{2q+1}\right)!}\right]^{-1} \sum_{s=0}^{2q} \left[A^{\frac{s}{2p+1}} \frac{t_{i+1}^{\frac{s-2q}{2q+1}}}{\frac{s-2q}{2q+1}!} - \sum_{k=0}^{\infty} \frac{(-1)^{k+1}}{k!} A^{\frac{s+(k+1)(2q+1)}{2p+1}} \cdot e^{\Delta_i A^{\frac{2q+1}{2p+1}}} \frac{\Delta_i^{k+\frac{s+1}{2q+1}}}{\frac{s-2q}{2q+1}!\left(k+\frac{s+1}{2q+1}\right)}\right] (19)$$

где $\Delta_i = t_{i+1} - t_i = const$.

Теперь остановимся на случае, когда начальная задача (6) заменяется задачей

$$D^\alpha x(t) = Ax(t) + f(t), \quad x(t_0) = x_0, \qquad (20)$$

где $f(t)$ на интервалах $(t_i, t_{i+1})$ являются постоянными, т.е. в интервале $(t_0, T)$ являются кусочно-постоянными. Тогда используя соотношения (19) из [9]

легко можно констатировать, что

$$x(t_{i+1}) = \Phi(t_{i+1}, t_i) x(t_i) + g(t_i), \qquad (21)$$

где $\Phi(t_{i+1}, t_i)$ определяется из (19), а $g(t_i)$ имеет вид [9]

$$g(t_i) = \left[ \sum_{s=0}^{2q} A^{\frac{s}{2p+1}} \frac{t_i^{\frac{s-2q}{2q+1}}}{\left(\frac{s-2q}{2q+1}\right)!} \right]^{-1} \sum_{s=0}^{2q} A^{\frac{s}{2p+1}} \left[ \left( \sum_{k=0}^{\infty} (-1)^k \frac{A^{\frac{s+(k+1)(2q+1)}{2p+1}}}{\frac{s-2q}{2q+1}! k! \left(k + \frac{s+1}{2q+1}\right)} \sum_{l=0}^{\infty} \frac{A^{\frac{l\,2q+1}{2p+1}}}{l!} \times \right.$$

$$\left. \times \frac{(t_{i+1} - t_i)^{l+k+1+\frac{s+1}{2q+1}}}{l+k+1+\frac{s+1}{2q+1}} \right) + \left( \frac{t_{i+1}^{\frac{s+1}{2q+1}} - t_i^{\frac{s+1}{2q+1}}}{\frac{s+1}{2q+1}!} \right) \right] f(t_i) \qquad (22)$$

Отметим, что формулу (21), (22) можно распространять для дискретизации линейных стационарных систем управления [11], которая в отличие [6] позволит точно решать соответствующие задачи управления. Специфика дискретизации классических задач и данной задачи являются существенно различными. Если в задаче (1) фундаментальная матрица $\Phi$ при $\Delta_i = \Delta$ является постоянной, то здесь, как видно из (19) и (22) фундаментальная матрица не является постоянной.

5. Численные сравнения

Теперь остановимся на сравнении результатов приближенных и точных фундаментальных матрицах (11) и (19). Поэтому анализируем их на следующем примере.

Пример 1. Рассмотрим следующее уравнение

$$D^{2\alpha} y + y = 0, \qquad 0 < \alpha < 1, \qquad (\text{п.1})$$

где с помощью соответствующих преобразований

$$y_1 = y, \qquad D^{\alpha} y_1 = y_2 \qquad (\text{п.2})$$

Получим системы дробно-производных дифференциальных уравнений

$$D^{\alpha} \begin{bmatrix} y_1 \\ y_2 \end{bmatrix} = \begin{pmatrix} 0 & 1 \\ -1 & 0 \end{pmatrix} \begin{bmatrix} y_1 \\ y_2 \end{bmatrix}. \qquad (\text{п.3})$$

Из за отсутствия нахождения фундаментальных матриц [12] пытаемся провести сравнение с приближенным результатом [6].

Пусть уравнение (п.3) определено в интервале [0.01, 1.01]. Разделим этот интервал на пять частей с постоянным шагом 0.2. Тогда по формуле (11), (12) имеем (при $\alpha = \dfrac{1}{3}$)

Фундаментальные матрицы по [6,(3.109)]

$$\bar{\Phi}(0.21,\ 0.01) = \begin{bmatrix} 3.333333333333333e-001 & 9.181368809065e-001 \\ -2.514702143092399e-001 & 3.333333333333e-001 \end{bmatrix},$$

$$\bar{\Phi}(0.41,\ 0.21) = \begin{bmatrix} -8.661856002099e-003 & 6.120376e-001 \\ -1.6764680953933e-001 & -8.6618529e-003 \end{bmatrix},$$

$$\bar{\Phi}(0.61,\ 0.41) = \begin{bmatrix} -5.804457201e-002 & 2.98092855e-001 \\ -8.16452e-002 & -5.8044572e-002 \end{bmatrix},$$
(п4)

$$\bar{\Phi}(0.81,\ 0.61) = \begin{bmatrix} -3.354369883e-002 & 1.70756678946e-001 \\ -4.67688637e-002 & -3.354369886e-002 \end{bmatrix},$$

$$\bar{\Phi}(1.01,\ 0.81) = \begin{bmatrix} -1.72097833e-002 & 1.3480937e-001 \\ -3.6923205e-002 & -1.72093e-002 \end{bmatrix}$$

Теперь находим точные фундаментальные матрицы по формуле (19) и имеем

$$\Phi(0.21,\ 0.01) = \begin{bmatrix} 8.253515142e-002 & 6.623231589e-002 \\ 6.6232315e-002 & 8.2535151422e-002 \end{bmatrix},$$

$$\Phi(0.41,\ 0.21) = \begin{bmatrix} 4.014333789e-001 & 3.91934853e-002 \\ -3.91934853771e-002 & 4.014333789e-001 \end{bmatrix},$$

$$\Phi(0.61,\ 0.41) = \begin{bmatrix} 3.9781124825e-001 & 5.4958790553e-002 \\ -5.495879055e-002 & 3.9781124811e-001 \end{bmatrix},$$
(п5).

$$\Phi(0.81,\ 0.61) = \begin{bmatrix} 3.957906410e-001 & 9.46203741871e-002 \\ -9.46203741873e-002 & 3.9579064106e-001 \end{bmatrix},$$

$$\Phi(1.01,\ 0.81) = \begin{bmatrix} 4.0040509424e-001 & 1.3019825994e-001 \\ -1.3019825995e-001 & 4.0040509421e-001 \end{bmatrix}$$

Вычислим следующие нормы между приближенными и точными фундаментальными матрицами в следующем виде

Составим таблицы фундаментальных матриц из (п4) с точностью $10^{-4}$

$$\bar{\Phi} = \begin{bmatrix} 1 & 0 & 0.3 & 0.9 & -0.0008 & 0.61 & -0.0058 & 0.029 & -0.033 & 0.017 & -0.017 & 0.134 \\ 0 & 1 & -0.25 & 0.3 & -0.167 & -0.0008 & -0.8016 & -0.0058 & -0.046 & 0.033 & -0.0369 & -0.017 \end{bmatrix}$$

Составим таблицы фундаментальных матриц из (п5) с точностью $10^{-4}$

$$\Phi = \begin{bmatrix} 1 & 0 & 0.08 & 0.06 & -0.0401 & 0.391 & 0.0397 & 0.0054 & 0.039 & 0.0946 & 0.0404 & 0.130 \\ 0 & 1 & -0.062 & 0.08 & -0.0397 & -0.0401 & -0.054 & 0.0397 & -0.094 & 0.039 & -0.13 & 0.0404 \end{bmatrix}$$

Разность фундаментальных матриц

По формуле (12) [6,(3.109)]

$$x(0) = \begin{bmatrix} 1 \\ 2 \end{bmatrix}; x(1) = \begin{bmatrix} 2.69 \\ 0.41 \end{bmatrix}; x(2) = \begin{bmatrix} 1.21 \\ -0.018 \end{bmatrix}; x(3) \begin{bmatrix} 0.53 \\ -0.019 \end{bmatrix}; x(4) = \begin{bmatrix} 0.30 \\ -0.011 \end{bmatrix}; \ x(5) = \begin{bmatrix} 0.25 \\ -0.0071 \end{bmatrix};$$

и по формуле (21)

$$y(0) = \begin{bmatrix} 1 \\ 2 \end{bmatrix}; y(1) = \begin{bmatrix} 0.21 \\ 0.098 \end{bmatrix}; y(2) = \begin{bmatrix} 0.091 \\ 0.031 \end{bmatrix}; y(3) \begin{bmatrix} 0.037 \\ 0.0079 \end{bmatrix}; y(4) = \begin{bmatrix} 0.01 \\ -0.0005 \end{bmatrix}; y(5) = \begin{bmatrix} 0.006 \\ -0.0022 \end{bmatrix};$$

Разность между решениями $\|x - y\| = 2.34$

Заключение

Приводиться новая формула определение фундаментальных матриц для систем линейных однородных дробно-производных дифференциальных уравнений с постоянными коэффициента. Далее результаты переноситься к случаю с возмущениями такие подходы позволяют дискретизировать задачи Коши.

**ЛИТЕРАТУРА**

# METHOD OF DISCRETIZING OF FRACTIONAL-DERIVATIVE LINEAR SYSTEMS OF ORDINARY DIFFERENTIAL EQUATIONS WITH CONSTANT COEFFICIENTS

**Aliev Fikret A., Aliev N. A., Velieva N.I., Gasimova K.G., Mamedova Y.V.**

**Abstract.** An exact discretization method is being developed for solving linear systems of ordinary fractional-derivative differential equations with constant matrix coefficients (LSOFDDECMC). It is shown that the obtained linear discrete system in this case does not have constant matrix coefficients. Further, this method is compared with the known approximate method. The above scheme is developed for arbitrary linear systems with piecewise constant perturbations. The results are applied to the discretization of linear controlled systems and are illustrated with numerical examples.